\documentclass{amsart}

\usepackage{amsfonts}
\usepackage{amsmath}
\usepackage{amsthm}
\usepackage{amssymb}
\usepackage[english]{babel}
\usepackage{graphics}
\usepackage{latexsym}
\usepackage{longtable} 
\usepackage{mathrsfs} 
\usepackage{graphicx}
\usepackage{wrapfig} 
\usepackage[pdftex,colorlinks=true,linkcolor=blue,citecolor=magenta]{hyperref}
\usepackage{enumitem}
\usepackage{esint}

\newtheorem{theorem}{Theorem}
\newtheorem{corollary}[theorem]{Corollary}
\newtheorem{lemma}[theorem]{Lemma}
\newtheorem{proposition}[theorem]{Proposition}

\newtheorem{claim}[theorem]{Claim}
\newtheorem{example}[theorem]{Example}

\theoremstyle{definition}
\newtheorem{definition}[theorem]{Definition}
\newtheorem{remark}[theorem]{Remark}

\newenvironment{dedication}
  {
   \itshape             
   \raggedleft          
  }
  {\par 
  \medskip \medskip
  }

\newcommand{\mL}{\mathcal{L}}
\newcommand{\mH}{\mathcal{H}}
\newcommand{\mF}{\mathcal{F}}

\newcommand{\mM}{\mathcal{M}}

\newcommand{\E}{\mathrm{E}}
\newcommand{\D}{\mathrm{D}}

\newcommand{\F}{\mathrm{F}}
\newcommand{\K}{\mathrm{K}}

\newcommand{\R}{\mathbb{R}}

\newcommand{\N}{\mathbb{N}}

\newcommand{\mB}{\mathbb{B}}

\renewcommand{\L}{\mathrm{L}}

\newcommand{\noi}{\noindent}
\newcommand{\ms}{\medskip}

\newcommand{\al}{\alpha}
\newcommand{\be}{\beta}
\newcommand{\ga}{\gamma}

\newcommand{\De}{\Delta}
\newcommand{\e}{\varepsilon}
\newcommand{\si}{\sigma}

\newcommand{\Om}{\Omega}

\newcommand{\av}{-\hspace{-10.5pt}\int}

\newcommand{\weak }{\, -\!\!\!\!-\!\!\!\!\rightharpoonup}
\newcommand{\weakstar }{ \overset{\, *_{\phantom{|}}}{{\smash{\weak }}\, } }

\newcommand{\larrow}{\longrightarrow}
\newcommand{\ot}{\otimes}

\newcommand{\LL}{\text{\LARGE$\llcorner$}}
\newcommand{\p}{\partial}
\newcommand{\sub}{\subseteq}
\newcommand{\set}{\setminus}
\newcommand{\by}{\times}

\newcommand{\ess}{\mathrm{ess}}

\renewcommand{\div}{\mathrm{div}}

\newcommand{\bt}{\begin{theorem}}\newcommand{\et}{\end{theorem}}
\newcommand{\bd}{\begin{definition}}\newcommand{\ed}{\end{definition}}
\newcommand{\bl}{\begin{lemma}}\newcommand{\el}{\end{lemma}}
\newcommand{\beq}{\begin{equation}}\newcommand{\eeq}{\end{equation}}
\newcommand{\bc}{\begin{claim}}\newcommand{\ec}{\end{claim}}
\newcommand{\bex}{\begin{example}}\newcommand{\eex}{\end{example}}
\newcommand{\bcor}{\begin{corollary}}\newcommand{\ecor}{\end{corollary}}
\newcommand{\bp}{\begin{proof}}\newcommand{\ep}{\end{proof}}

\numberwithin{equation}{section}

\begin{document}

\title[$L^\infty$ fully nonlinear inverse source identification]{ON THE INVERSE SOURCE IDENTIFICATION PROBLEM IN $L^\infty$ FOR FULLY NONLINEAR ELLIPTIC PDE}
 
 \author{Birzhan Ayanbayev}

\address{Freie Universit\"at Berlin, Arnimallee 6, 14195 Berlin, Germany}

\email{bayanbayev@zedat.fu-berlin.de}

\author{Nikos Katzourakis}


\address{Department of Mathematics and Statistics, University of Reading, Whiteknights, PO Box 220, Reading RG6 6AX, United Kingdom}

\email{(Corresponding author) n.katzourakis@reading.ac.uk}

\begin{dedication}
\small Dedicated to Enrique Zuazua on the occasion of his 60th birthday
\end{dedication}
  


\begin{abstract} In this paper we generalise the results proved in [N. Katzourakis, \emph{An $L^\infty$ regularisation strategy to the inverse source identification problem for elliptic equations}, SIAM J. Math. Anal. 51:2, 1349-1370 (2019)] by studying the ill-posed problem of identifying the source of a fully nonlinear elliptic equation. We assume Dirichlet data and some partial noisy information for the solution on a compact set through a fully nonlinear observation operator. We deal with the highly nonlinear nonconvex nature of the problem and the lack of weak continuity by introducing a two-parameter Tykhonov regularisation with a higher order $L^2$ ``viscosity term" for the $L^\infty$ minimisation problem which allows to approximate by weakly lower semicontinuous cost functionals.
\end{abstract}

\maketitle


\section{Introduction} \label{section1}

Let $n,k\in\N$ with $k,n\geq 2$ and let $\Om \sub \R^n$ be a bounded connected domain with $C^{1,1}$ regular boundary $\p\Om$. Let also
\[
F\ : \ \ \Om \by \R \by \R^n \by \R^{n^{\ot 2}}_s \larrow \R
\]
be a Carath\'eodory function, namely $x\mapsto F(x,\mathrm r,\mathrm{p}, \mathrm X)$ is Lebesgue measurable for all $(\mathrm r, \mathrm p, \mathrm X) \in \smash{\R \by \R^n \by \R^{n^{\ot 2}}_s}$ and $(\mathrm r, \mathrm p, \mathrm X)\mapsto F(\mathrm x, \mathrm r, \mathrm p, \mathrm X)$ is continuous for a.e.\ $x \in\Om$. In this paper, the notation $\R^{n^{\ot k}}_s$ stands for the vector space of fully symmetric $k$-th order tensors in $\R^n \ot \cdots \ot \R^n$ ($k$-times). Given $g\in W^{2,\infty}(\Om)$, consider the Dirichlet problem
\beq  \label{1.3}
\left\{
\begin{array}{rl}
\F [u]\, =\, f, & \text{ in }\Om,
\\
u\, =\, g, & \text{ on }\p \Om,
\end{array}
\right.
\eeq
for some appropriate source $f :\Om\larrow \R$. Here $\F[u]$ denotes the induced fully nonlinear 2nd order differential operator, defined on smooth functions $u$ as
\beq \label{1.1}
\F [u]\,:=\, F(\cdot, u, \D u, \D^2 u).
\eeq
Evidently, we are employing the standard symbolisations $\D u =(\D_i u)_{i=1...n}$, $\D^2 u =(\D^2_{ij} u)_{i,j=1...n}$ and $\D_i\equiv \p/\p x_i$. The above direct Dirichlet problem for $\F$ asks to determine $u$, given a source $f$ and boundary data $g$. (In fact the source $f$ is obsolete and can be absorbed into $F$, but for the problem we are interested in this paper it is more convenient to write it in this separated form). This is a semi-classical problem which is essentially standard material, see e.g.\ \cite{GT}. In particular, it is known that under various sets of assumptions on $F$ that \eqref{1.3} is well-posed and, given $f\in L^\infty(\Om)$ and $g\in W^{2,\infty}(\Om)$, for any $p>n$ there exists a unique solution $u$ in the locally convex (Fr\'echet) space
\[
\mathcal{W}^{2,\infty}_g(\Om)\,:= \bigcap_{1<p<\infty} \big(W^{2,p}\cap W^{1,p}_g\big)(\Om)  .
\]
In general, the solution $u$ is not in the smaller space $\smash{(W^{2, \infty}\cap W^{1,\infty}_g)(\Omega)}$ (not even locally), due to the failure of the $W^{2,p}$ estimates for $p=\infty$, which happens even in the linear case (see e.g.\ \cite{GM}). Additionally, \eqref{1.3} satisfies for any $p>n$ the fully nonlinear $L^p$ global estimate
\beq \label{estimate}
\|F(\cdot,v, \D v,\D^2 v)\|_{L^p(\Om)}\ge C_1\|v\|_{W^{2,p}(\Om)} -C_1\|g\|_{W^{2,p}(\Om)} -C_2 
\eeq
for some constants $C_1,C_2>0$ depending only on the parameters and any $v\in \big(W^{2,p}\cap W^{1,p}_g\big)(\Om)$. For sufficient conditions on $F$ which guarantee the satisfaction of solvability of \eqref{1.2} in the strong sense and of the uniform estimate \eqref{estimate} we refer to \cite{Caffarelli, Caffarelli-Cabre,  Dong-Krylov-Li, Krylov, Li-Zhang}.

Note that the above problem contains as a special case the archetypal instance of divergence operators with $C^1$ matrix coefficient $A$, as well as the non-divergence linear case with continuous coefficient:
\beq \label{linear-pdes}
\left\{
\begin{split}
\L_1 [u]\, & =\, \div(A \D u)\,+\, b \cdot \D u \,+\, cu,
\\
\L_2 [u]\, & =\, A :\D^2 u\,+\, b \cdot \D u \,+\, cu.
\end{split}
\right.
\eeq
In the above, the notations ``$:$" and ``$\cdot$" symbolise the Euclidean inner products in the space of symmetric matrices $\R^{n^{\ot 2}}_s$ and in $\R^n$ respectively. More generally, the inner product of two tensors $T,S \in \R^{n^{\ot k}}_s$ will also be denoted by ``$:$", that is
\[
T:S\, := \sum_{1\leq a_1,\ldots ,a_k \leq n} T_{a_1\cdots a_k} \, S_{a_1\cdots a_k}.
\]

The {\it inverse problem} relating to \eqref{1.3} asks the question of perhaps determining $f$, given the boundary data $g$ and some other {\it partial information on the solution $u$, typically some approximate experimental measurements of some function of it known only up to some error.} The inverse problem is {\it severely ill-posed} even in the linear case of the Laplacian operator $\F=\De$, as the noisy data measured on a subset of $\Om$ might either not be compatible with {\it any} exact solution, or they may not suffice to determine a unique $f$ even if compatibility holds true.

The above type of inverse problems are especially crucial for various applications, even in the model case of the Poisson equation, see e.g.\ \cite{AAM, BD, EHN, I, LHY, MV, NA, SH, X, Y, YF, ZMYX}. In this paper we will assume that the approximate information on $u$ takes the form
\beq
\label{1.5}
\K[u]\, =\, k^\ga \ \ \text{ on }\mathcal{K},
\eeq
where $\K$ is an {\it observation operator}, taken to be a first order fully nonlinear differential operator of the form
\beq
\label{1.6}
\K[u]\,:=\, K(\cdot,u,\D u),
\eeq
where $K$ and its partial derivates $K_{\mathrm r}, K_{\mathrm p}$ satisfy
\beq
\label{1.7}
K,\, K_{\mathrm r}\in C(\mathcal{K} \by \R \by \R^n),
\ \ \ \
 K_{\mathrm p} \in C\big(\mathcal{K} \by \R \by \R^n;\R^n\big).
\eeq
In \eqref{1.5} and \eqref{1.7}, $\mathcal{K}$ symbolises the set on which we take measurements, which will be assumed to satisfy
\beq
\label{1.8}
\mathcal{K} \sub \overline{\Om} \text{ is compact and }\text{exists }\kappa\in[0,n]: \ \mH^\kappa(\mathcal{K})<\infty.
\eeq
In the above, $\mH^\kappa$ denotes the Hausdorff measure of dimension $\kappa$. Our general measure and functional notation will be either standard or self-explanatory, e.g.\ as in \cite{D,E,KV}. Finally, $k^\ga \in L^\infty(\mathcal{K},\mH^\kappa)$ is the function of approximate (deterministic) measurements taken on $\mathcal{K}$, at noise level at most $\ga>0$:
\beq 
\label{1.9}
\| k^\ga - k^0 \|_{L^\infty(\mathcal{K},\mH^\kappa)} \, \leq\, \ga,
\eeq
where $k^0=\K[u^0]$ corresponds to ideal error-free measurements of an exact solution to \eqref{1.3} with source $f=\F[u^0]$. 

To recapitulate, in this paper we study the following ill-posed inverse source identification problem for fully nonlinear elliptic PDEs: 
\beq
\label{1.10}
\left\{
\begin{array}{rl}
\F [u]\, =\, f\ , & \text{ in }\Om,
\\
u\, =\, g\ , & \text{ on }\p \Om,
\\
\K[u]\, =\, k^\ga, & \text{ on }\mathcal{K}.
\end{array}
\right.
\eeq
This means that we are searching for a selection process of a suitable approximation for $f$ from the data $k^\ga$ on $\mathcal{K}$ through the observation $\K[u]$ of the solution $u$. To the best of our knowledge, \eqref{1.10} has not been studied before, at least in this generality. Our approach does not exclude the extreme cases of $\mathcal{K}=\overline{\Om}$ (full information) or of $\mathcal{K}=\emptyset$ (no information), although trivial changes are required in the proofs. Sadly, an exact solution may not exist as the constraint may be incompatible with the solution of \eqref{1.3}, owing to the errors in measurements. On the other hand, it is not possible to have a uniquely determined source on the constraint-free region $\Om\set\mathcal{K}$.  Instead, our goal is a strategy to determine an optimally fitting $u^\ga$ (and respective source $f^\ga:=\F[u^\ga]$) to the ill-posed problem \eqref{1.10}. A popular choice of operator $\K$ in the literature (when $\L=\De$) consist of some term of the separation of variables formula, as e.g.\ in \cite{YF}. 

Herein we follow an approach based on recent advances in Calculus of Variations in the space $L^\infty$ (see e.g.\ \cite{KM} and references therein) developed for functionals involving higher order derivatives, which has already been applied to the special case of the inverse source problem for linear PDEs \eqref{linear-pdes} in \cite{K2}. This relatively new field field was pioneered by Gunnar Aronsson in the 1960s  (see e.g.\ \cite{A1,A2,A3,A,AB}) and is still a very active area of research; for a review of the by-now classical theory involving scalar first order functionals we refer to \cite{K1}. 

Following \cite{K2}, we aim at providing a {\it regularisation strategy} inspired by the classical Tykhonov regularisation strategy in $L^2$ (see e.g.\ \cite{Ki,N}). As a first possible step, consider the next putative $L^\infty$ ``error" functional:
\beq 
 \label{1.11-0} 
\ \ \mathrm{E}_{\infty,\al} (u)\, :=\,\big \|\K[u]-k^\ga \big\|_{L^\infty(\mathcal{K},\mH^\kappa)}+\, \al \big \| \F[u] \big\|_{L^\infty( \Om)}, \ \ u\in  \mathcal{W}^{2,\infty}_g (\Om),
\eeq
for some fixed parameter $\al>0$. The advantage of searching for a best fitting solution in $L^\infty$ is evident: we can keep the error term $|\K[u]-k^\ga |$ uniformly small and not merely small on average, as would happen if one chose to minimise some integral of a power of the error instead. As in \cite{K2}, the goal would be to minimise $\mathrm{E}_{\infty,\al}$ over $\mathcal{W}^{2,\infty}_g (\Om)$, and then any minimiser of \eqref{1.11-0} would provide a candidate solution for our problem. Then, for any fixed $\al$, this would be the best fitting solution with the least possible uniform error, namely $\F[u]\cong f$ uniformly on $\Om$ and $\K[u]\cong k^\ga$ uniformly on $\mathcal{K}$. Unfortunately, even if one momentarily ignores the problem of lack of regularity for \eqref{1.11-0} and the fact that $\mathcal{W}^{2,\infty}_g (\Om)$ is not a Banach space, the main problem is that in general minimisers do not exist in the genuine fully nonlinear case of operator $\F$ (namely when $\mathrm X \mapsto \F(\mathrm x, \mathrm r, \mathrm p, \mathrm X)$  is nonlinear) as \eqref{1.11-0} is not weakly lower semicontinuous in the Fr\'echet space $\mathcal{W}^{2,\infty}_g (\Om)$, as the highest order term may be nonconvex/non-quasiconvex. In the special linear case of \cite{K2}, this problem was not present as the linearity of the differential operator was implying the desired weak lower semi-continuity.

In this work we resolve the problem explained above by proposing a double approximation method (or rather triple, as we will see shortly) which involves an additional Tykhonov or ``viscosity" term which effectively is a weakly lower-semicontinuous approximation of \eqref{1.11-0}. Hence, we will consider instead
\beq 
 \label{1.11} 
\left\{
\begin{split}
\mathrm{E}_{\infty,\al,\be} (u)\,  :=\,\big \|  &  \K[u]-k^\ga \big\|_{L^\infty(\mathcal{K},\mH^\kappa)}  +\, \al \big \| \F[u] \big\|_{L^\infty( \Om)}+\,\frac{\be}{2} \|\D^{\bar n} u\|^2_{L^2(\Om)}, \phantom{\bigg|} \!\!\!\!
\\
& \ \ \ \ \ \ u\in \big(W^{\bar n,2} \cap W^{1,2}_g\big)(\Om).
\end{split}
\right.
\eeq
where $\bar n:= [n/2]+3$. In the above $\be>0$ is a fixed parameter, $[\,\cdot\,]$ symbolises the integer part and $\D^{\bar n} u$ is the $\bar n$-th order weak derivative of $u$. 

It is well known in the Calculus of Variations in $L^\infty$ that (global) minimisers of supremal functionals, although usually simple to obtain with a standard direct minimisation (\cite{D,FL}), they are not genuinely minimal as they do not share the nice ``local" optimisation properties of minimisers of their integral counterparts (see e.g.\ \cite{BN, RZ}). The case of \eqref{1.11} studied herein is no exception to this rule. A relatively standard method is bypass these obstructions is to employ minimisers of $L^p$ approximating functionals as $p\to \infty$, establishing appropriate convergence of $L^p$ minimisers to a limit $L^\infty$ minimiser.  The idea underlying this approximation technique is based on the simple measure theory fact that the $L^p$ norm (of a function in $L^1\cap L^\infty$) converges to the $L^\infty$ norm as $p\to \infty$.  This method is quite standard in the field and {\it furnishes a selection principle of $L^\infty$ minimisers with additional desirable properties} (see e.g.\  \cite{BJW1, BP, CDP, GNP, KM}). In this fashion one is also able to bypass the lack of differentiability of supremal functionals and derive necessary PDE conditions satisfied by $L^\infty$ extrema. This is indeed the method that is employed in this work as well, along the lines of \cite{K2}.

We now present the main results to be established in this paper. As already explained, we will obtain {\it special} minimisers of \eqref{1.11} as limits of minimisers of
 \beq 
 \label{1.12} 
 \left\{
\begin{split}
\mathrm{E}_{p,\al,\be} (u) := \big \|   |\K[u]-  &  k^\ga|_{(p)}  \big\|_{L^p(\mathcal{K},\mH^\kappa)}+ \al \big \| |\F[u]|_{(p)}  \big\|_{L^p( \Om)}+ \frac{\be}{2} \|\D^{\bar n} u\|^2_{L^2(\Om)},  \phantom{\bigg|}
\\
& \ \ u\in \big(W^{\bar n,2} \cap W^{1,2}_g\big)(\Om)  .
\end{split}
\right. \!\!\!\!\!\!\!\!
\eeq
In \eqref{1.12} we have used the normalised $L^p$ norms 
\[
\big \| f\big\|_{L^p(\mathcal{K},\mH^\kappa)}\, :=\, \left(\, {\, {\av_\mathcal{K}}} |f|^p \, \mathrm{d}\mH^\kappa\right)^{1/p}, \ \ \ \big \| f\big\|_{L^p(\Om)}\, :=\, \left(\, {\, {\av_\Om}} |f|^p \, \mathrm{d}\mL^n\right)^{1/p}
\]
and the integral signs with slashes symbolises the average with respect to the Hausdorff measure $\mH^\kappa$ and the Lebesgue measure $\mL^n$ respectively. Further, $|\, \cdot \,|_{(p)}$ symbolises the next $p$-regularisation of the absolute value away from zero:
\[
|a|_{(p)}\, :=\, \sqrt{|a|^2+p^{-2}}.
\]
Let us also note that, due to our $L^p$-approximation method, as an auxiliary result we also provide an $L^p$ regularisation strategy for finite $p$ as well, which is of independent interest. For the proof we will need to assume that the Dirichlet problem \eqref{1.3} for the fully nonlinear operator $\F$ satisfies the $W^{2,p}$ elliptic estimates \eqref{estimate} for all large enough (finite) $p>n$, as well the following: $F$ and its partial derivates $F_{\mathrm r}, F_{\mathrm p}, F_{\mathrm X}$ satisfy
\beq \label{1.2}
\left\{\ \ 
\begin{split}
& F,\, F_{\mathrm r} \in  C\big(\overline{\Om} \by \R \by \R^n \by \R^{n^{\ot 2}}_s\big),
\\
& F_{\mathrm r} \in   C\big(\overline{\Om} \by \R \by \R^n \by \R^{n^{\ot 2}}_s,\R^n\big),
\\
& F_{\mathrm X} \in  C\big(\overline{\Om} \by \R \by \R^n \by \R^{n^{\ot 2}}_s ; \R^{n^{\ot 2}}_s\big).
\end{split}
\right.
\eeq
We note that sufficient general conditions of when such operators satisfy $W^{2,p}$ elliptic estimates can be found for instance in the papers \cite{Caffarelli, Caffarelli-Cabre,  Dong-Krylov-Li, Krylov, Li-Zhang}.

\bt
\label{theorem1}
Let $\Om \sub \R^n$ be a bounded $C^{1,1}$ domain and $g \in W^{\bar n, 2}(\Om)$. Suppose also the operators \eqref{1.1} and \eqref{1.6} are given, satisfying the assumptions \eqref{estimate}, \eqref{1.7}, \eqref{1.8} and \eqref{1.2}. Suppose further a function $k^\ga \in L^\infty(\mathcal{K},\mH^\kappa)$ is given which satisfies \eqref{1.9} for $\ga>0$. Let finally $\al,\be>0$ be fixed. Then, we have the following results for the inverse problem associated to \eqref{1.10}:

\smallskip

\noi {\bf (i)} There exists a global minimiser $u_\infty\equiv u_\infty^{\al,\be,\ga} \in \big(W^{\bar n,2} \cap W^{1,2}_g\big)(\Om)$ of the functional $\E_{\infty,\al,\be}$ defined in \eqref{1.11}. In particular, we have $E_\infty(u_\infty)\leq E_\infty (v)$ for all $v\in \big(W^{\bar n,2} \cap W^{1,2}_g\big)(\Om)$ and 
\[
f_\infty \equiv f^{\al,\be,\ga}_\infty \,:= \,\F[u_\infty^{\al,\be,\ga}] \, \in L^\infty(\Om).
\]
Further, there exist signed Radon measures 
\[
\mu_\infty \equiv \mu_\infty^{\al,\be,\ga} \, \in \, \mM(\Om), \  \  \ \nu_\infty \equiv \nu_\infty^{\al,\be,\ga} \, \in \, \mM(\mathcal{K})
\]
such that the nonlinear divergence PDE
\beq
\label{1.14}
K_{\mathrm r} [u_\infty]\nu_\infty \, -\, \div\big( \K_{\mathrm p} [u_\infty]\nu_\infty  \big)\, +\, \al ((\mathrm{d} \F)_{u_\infty})^*[\mu_\infty]\,+\,\be (-1)^{\bar n} (\D^{\bar n} : \D^{\bar n} u_\infty)   =\, 0, \phantom{\Big|}
\eeq
is satisfied by the triplet $(u_\infty,\mu_\infty,\nu_\infty)$ in the sense of distributions (see \eqref{1.16-A}). In \eqref{1.14}, the operator $((\mathrm{d} \F)_{u_\infty})^*$ is the formal adjoint of the linearisation of $\F$ at $u_\infty$, defined via duality as
\[
\begin{split}
((\mathrm{d} \F)_{u_\infty})^*[v]\,:=\, \div(\div (\F_{\mathrm X} [u_\infty] v))\, -\, \div(\F_{\mathrm p} [u_\infty] v)\,+\,\F_{\mathrm r} [u_\infty]v,
\end{split}
\]
$F_{\mathrm r} ,K_{\mathrm r} ,F_p,K_p,F_X$ denote the partial derivatives of $F,K$ with respect to the respective variables and $\F_{\mathrm r} [v],\K_{\mathrm r} [v],\F_{\mathrm p} [v],\K_{\mathrm p} [v],\F_{\mathrm X} [v]$ denote the respective differential operators $F_{\mathrm r} (\cdot,v,\D v,\D^2 v), K_{\mathrm r} (\cdot,v,\D v),...$ etc. Additionally, the error measure $\nu_\infty$ is supported in the closure of the subset of $\mathcal{K}$ of maximum error, namely
\beq
\label{1.15}
\mathrm{supp} (\nu_\infty) \,\sub \, \Big\{  \big|\K[u_\infty]-k^\ga \big|^\bigstar = \big\| \K[u_\infty]-k^\ga \big\|_{L^\infty(\mathcal{K},\mH^\kappa)} \Big\},
\eeq
where  $``\, (\, \cdot\, )^\bigstar \,"$ denotes the ``essential limsup" with respect to $\mH^\kappa\LL_\mathcal{K}$ on $\mathcal{K}$ (see Definition \ref{def} that follows). If additionally the data function $k^\ga$ is continuous on $\mathcal{K}$, \eqref{1.15} can be improved to
\beq
\label{1.16}
\mathrm{supp} (\nu_\infty) \,\sub \, \Big\{  \big|\K[u_\infty]-k^\ga \big| = \big\| \K[u_\infty]-k^\ga \big\|_{L^\infty(\mathcal{K},\mH^\kappa)} \Big\} .
\eeq

\ms

\noi {\bf (ii)} For any $\al,\be,\ga>0$, the minimiser $u_\infty$ can be approximated by a family of minimisers $(u_p)_{p>n}\equiv(u_p^{\al,\be,\ga})_{p>n}$ of the respective $L^p$ functionals \eqref{1.12} and the pair of measures $(\mu_\infty,\nu_\infty)  \in \mM(\Om)\by \mM(\mathcal{K})$ can be approximated by respective signed measures $(\mu_p,\nu_p)_{p>n}\equiv (\mu^{\al,\be,\ga}_p,\nu^{\al,\be,\ga}_p)_{p>n}$, as follows: 

\noi For any $p>n$, \eqref{1.12} has a global minimiser $u_p\equiv u_p^{\al,\be,\ga}$ in $\big(W^{\bar n,2} \cap W^{1,2}_g\big)(\Om)$ and there exists a  sequence $p_j\larrow \infty$ as $j\to \infty$, such that
\beq
\label{1.17}
\left\{\ \ 
\begin{array}{ll}
u_p^m \larrow u_p, & \text{ in }C^{2}(\overline{\Om}),
\smallskip
\\
\D^k u_p^m \larrow \D^k u_p, & \text{ in }L^{2}\big(\Om,\R^{n^{\ot k}}_s\big), \text{ for all }k\in\{3,...,\bar n-1\},
\smallskip
\\
\D^{\bar n} u_p^m \weak \D^{\bar n} u_p, & \text{ in }L^{2}\big(\Om,\R^{n^{\ot \bar n}}_s\big),
\smallskip
\end{array}
\right.
\eeq
as $p_j\to \infty$. Moreover, we have
\beq
\label{1.18}
\left\{\ \ \
\begin{split}
\nu_p\, & :=\, \frac{\big|\K[u_p]-k^\ga \big|^{p-2}_{(p)} \big(\K[u_p]-k^\ga \big) }{ \mH^\kappa(\mathcal{K}) \, \big \| |\K[u_p]-k^\ga|_{(p)} \big\|^{p-1}_{L^p(\mathcal{K},\mH^\kappa)}}\mH^\kappa \LL_\mathcal{K} \, \weakstar \, \nu_\infty,  \ \ \ \text{ in }\mM(\mathcal{K}),
\\
\mu_p\, & :=\, \frac{| \F[u_p] |^{p-2}_{(p)} \,  \F[u_p] }{ \mL^n(\Om) \, \big \| | \F[u_p] |_{(p)} \big\|^{p-1}_{L^p(\Om)}}\mL^n \LL_\Om \, \weakstar \, \mu_\infty,  \hspace{50pt} \text{ in }\mM(\Om),
\end{split}
\right.
\eeq
as $p_j\to \infty$. Further, for each $p>n$,  the triplet $(u_p,\mu_p,\nu_p)$ solves the PDE
\beq
\label{1.19}
K_{\mathrm r} [u_p]\nu_p \, -\, \div\big( \K_{\mathrm p} [u_p]\nu_p  \big)\, +\, \al ((\mathrm{d} \F)_{u_p})^*[\mu_p]\,+\,\be  (-1)^{\bar n} (\D^{\bar n} : \D^{\bar n} u_p)   =\, 0, \phantom{\Big|}
\eeq
in the sense of distributions (see \eqref{1.21-A}).

\ms

\noi {\bf (iii)} For any exact solution $u^0 \in \big(W^{\bar n,2} \cap W^{1,2}_g\big)(\Om)$ of \eqref{1.10} (with $f=\F[u^0]$ and $\K[u^0]=k^0$) corresponding to measurements with zero error, we have the estimate: 
\beq
\label{1.20}
\Big\| \K[u^{\al,\be,\ga}_\infty] - \K[u^0] \Big\|_{L^\infty(\mathcal{K},\mH^\kappa)}  \leq\, 2\ga + \al \, \| \F[u^0]\|_{L^\infty(\Om)}\,+\, \frac{\be}{2} \, \|\D^{\bar n} u^0\|^2_{L^2(\Om)},
\eeq
for any $\al,\be,\ga>0$.

\ms

\noi {\bf (iv) }  For any exact solution $u^0 \in \big(W^{\bar n,2} \cap W^{1,2}_g\big)(\Om)$ of \eqref{1.10} (with $f=\F[u^0]$ and $\F[u^0]=k^0$) corresponding to measurements with zero error and for $p>n$, we have the estimate: 
\beq
\label{1.24}
\Big\| \K[u^{\al,\be,\ga}_p] - \K[u^0] \Big\|_{L^p(\mathcal{K},\mH^\kappa)}  \leq\, 2\ga + \al \, \| \F[u^0]\|_{L^p(\Om)}+\, \frac{\be}{2} \, \|\D^{\bar n} u^0\|^2_{L^2(\Om)},
\eeq
for any $\al,\be,\ga>0$.
\et 

We note that the estimate in part {\bf (iv)} is useful if we have merely that $\F[u^0] \in L^p (\Om)$ for $p<\infty$ (namely when $\F[u^0] \not\in L^\infty(\Om)$). We close this section by noting that the reader may find in \cite{K2} various comments and counter-examples regarding the optimality of Theorem \ref{theorem1} (therein stated for the case of a linear differential operator $\F$).

\ms

\section{Discussion, auxiliary results and proofs}

We begin with some clarifications on Theorem \ref{theorem1}.

\begin{remark} \label{remark4} In index form, the definition of the formal adjoint $((\mathrm{d} \F)_{u_\infty})^*$ of the linearisation of $\F$ at $u_\infty$ can be written as
\beq 
\label{dual-op}
(\mathrm{d} \F)_{u_\infty})^*[v]\,=\, \sum_{i,j=1}^n\D^2_{ij}(\F_{\mathrm X_{ij}}[u_\infty]v)\, -\, \sum_{k=1}^n\D_k(\F_{\mathrm p_{k}}[u_\infty]\, v)\,+\,\F_{\mathrm r}[u_\infty] \, v
\eeq
and its distributional interpretation via duality reads
\beq
\big\langle (\mathrm{d} \F)_{u_\infty})^*[v],\phi \big\rangle = \int_\Om \bigg\{\sum_{i,j=1}^n(\D^2_{ij}\phi) \F_{\mathrm X_{ij}}[u_\infty]\, +\, \sum_{k=1}^n(\D_k\phi) \F_{\mathrm p_{k}}[u_\infty]\,+\,\phi \F_{\mathrm r}[u_\infty] \bigg\} v \, \mathrm{d} \mL^n,
\eeq
for all $\phi \in C^{\bar n}_c(\Om)$. Hence, by taking into account the definitions of the measures $\mu_p,\nu_p$ in \eqref{1.18},  the distributional interpretation of \eqref{1.19} is
\beq \label{1.21-A}
\begin{split}
&\, {\av_\mathcal{K}} \Big( \K_{\mathrm r} [u_p] \,\phi \, +\, \K_{\mathrm p}[u_p]\cdot \D \phi \Big) \frac{\big|\K[u_p]-k^\ga \big|^{p-2}_{(p)} \big(\K[u_p]-k^\ga \big)}{\big \| |\K[u_p]-k^\ga|_{(p)} \big\|^{p-1}_{L^p(\mathcal{K},\mH^\kappa)}} \, \mathrm{d} \mH^\kappa \, 
\\
&+\ \al \, {\, {\av_\Om}} \Big(\F_{\mathrm r} [u_p]\phi+\F_{\mathrm p} [u_p]\cdot\D \phi+\F_{\mathrm X} [u_p]:\D^2 \phi\Big) \frac{| \F[u_p] |^{p-2}_{(p)} \,  \F[u_p]  }{ \big \| | \F[u_p] |_{(p)} \big\|^{p-1}_{L^p(\Om)}}  \, \mathrm{d} \mL^n\,
\\
&+\ \beta \, {\, {\av_\Om}} \, \D^{\bar n} u_p : \D^{\bar n} \phi \, \mathrm{d} \mL^n\, =\, 0,
\end{split}
\eeq
for all $\phi \in C^{\bar n}_c(\Om)$. Similarly,  the distributional interpretation of \eqref{1.14} is
\beq \label{1.16-A}
\begin{split}
&\, {\av_\mathcal{K}} \Big( \K_{\mathrm r} [u_\infty] \,\phi \, +\, \K_{\mathrm p} [u_\infty]\cdot \D \phi \Big) \, \mathrm{d} \nu_\infty 
\\
&+\ \al \, {\, {\av_\Om}} \Big(\F_{\mathrm r} [u_\infty]\phi+\F_{\mathrm p} [u_\infty]\cdot\D \phi+\F_{\mathrm X} [u_\infty]:\D^2 \phi\Big)  \, \mathrm d \mu_\infty\,
\\
&+\ \beta \, {\, {\av_\Om}} \, \D^{\bar n} u_\infty : \D^{\bar n} \phi \, \mathrm{d} \mL^n\, =\, 0,
\end{split}
\eeq
for all $\phi \in C^{\bar n}_c(\Om)$. 
\end{remark}

We now state a definition and a result taken from \cite{K2} which are required for our proofs. 

\begin{definition}[The essential limsup, \cite{K2}] \label{def} \ 
Let $X \sub \R^n$ be a Borel set and let $\nu \in \mM(X)$ be a finite positive Radon measure on $X$. Given $f\in L^\infty(X,\nu)$, we define $f^\bigstar \in L^\infty(X,\nu)$ by
\[
f^\bigstar(x) \, := \, \lim_{\e\to 0} \bigg(\nu-\underset{y \in \mB_\e(x)}{\ess\,\sup} \, f(y)\bigg), \ \ \ \forall\, x\in X,
\]
and call $f^\star$ {\bf the $\nu$-essential limsup of $f$}. Here $\mB_\e(x)$ denotes the open ball of radius $\e$ centred at $x\in X$. \end{definition}

The following result studies what we call ``concentration measures" of the approximate $L^p$ minimisation problems as $p\to \infty$.

\begin{proposition}[$L^p$ concentration measures as $p\to\infty$] \label{prop} Let $X$ be a compact metric space, endowed with a positive finite Borel measure $\nu$ which gives positive values to any open subset of $X$ except $\emptyset$. Consider $(f_p)_1^\infty \sub L^\infty(X,\nu)$ and the sequence of signed Radon measures $(\nu_p)_1^\infty \sub \mM(X)$, given by:
\[
\nu_p\,:=\, \frac{1}{ \nu(X)} \frac{\big(|f_p|_{(p)}\big)^{p-2} f_p }{ \big\| |f_p|_{(p)} \big\|^{p-1}_{L^p(X,\nu)}}\, \nu, \ \ \ p\in \N,
\]
where $|\cdot|_{(p)}=(|\cdot|^2+p^{-2})^{1/2}$. Then:
\ms

\noi {\bf (i)} There exists a subsequence $(p_i)_1^\infty$ and a limit measure $\nu_\infty \in \mM(X)$ such that
\[
\nu_p \, \weakstar \, \nu_\infty \ \ \text{ in }\mM(X),
\]
as $p_i\to \infty$.

\ms

\noi {\bf (ii)} If there exists $f_\infty \in L^\infty(X,\nu)\set \{0\}$ such that
\[
\sup_X |f_p -f_\infty| \larrow 0\ \ \text{ as }p\to\infty,
\]
then $\nu_\infty$ is supported in the set where $|f_\infty|$ is maximised:
\[
\mathrm{supp}(\nu_\infty) \, \sub \, \Big\{ |f_\infty|^\bigstar= \|f_\infty\|_{L^\infty(X,\nu)} \Big\}.
\]

\noi {\bf (iii)} If additionally to {\bf(ii)} the modulus $|f_\infty|$ of $f_\infty$ is continuous on $X$, then the following stronger assertion holds:
\[
\mathrm{supp}(\nu_\infty) \, \sub\, \Big\{ |f_\infty|= \|f_\infty\|_{L^\infty(X,\nu)} \Big\}.
\]
\end{proposition}

Now we establish Theorem \ref{theorem1}. The proof consists of several lemmas. We note that some details might be quite well known to the experts of Calculus of Variations, but we chose to give most of the niceties for the convenience of the readers and for the sake of completeness of the exposition.

\begin{lemma} \label{lemma1} For any $p>n$ and $\al,\be,\ga>0$, the functional \eqref{1.12} has a minimiser $u_p \in (W^{\bar n, 2}\cap W^{1,2}_g )(\Om)$:
\[
\E_{p,\al,\be}(u_p) \, =\, \inf\Big\{\E_{p,\al,\be}(v)\ : \ v\in (W^{\bar n, 2}\cap W^{1,2}_g )(\Om) \Big\}.
\]
\end{lemma}

\bp Since $g\in W^{2,\infty}(\Om)$ 
\[
\begin{split}
\E_{p,\al,\be}(g)\, \leq \,&\, \E_{\infty,\al,\be}(g)\\
 \leq \, & \, \|k^\ga \|_{L^\infty(\mathcal{K},\mH^\kappa)}\, +\  \|\K(\cdot,g,\D g) \|_{L^\infty(\mathcal{K},\mH^\kappa)}
\\
& +\al \|\F(\cdot,g, \D g,\D^2 g)\|_{L^\infty(\Om)} + \frac{\beta}{2} \|\D^{\bar n} g\|^2_{L^2(\Om)}\\
<\, &\, \infty.
\end{split}
\]
Hence, 
\[
0\, \leq\, \inf\Big\{\E_{p,\al,\be}(v)\ : \ v\in (W^{\bar n, 2}\cap W^{1,2}_g )(\Om) \Big\}\, \leq\, \E_{\infty,\al,\be}(g) \,<\,\infty.
\]
Further, $\E_{p,\al,\be}$ is coercive in the space $(W^{\bar n, 2}\cap W^{1,2}_g )(\Om)$. Indeed, by our assumption \eqref{estimate}, H\"older's inequality and that $p>n\geq 2$, for any $v\in (W^{\bar n, 2}\cap W^{1,2}_g )(\Om)$ we have
\[
\begin{split}
\E_{p,\al,\be}(v) \, &\geq\, \al \| \F[v] \|_{L^p(\Om)} + \frac{\beta}{2} \|\D^{\bar n} v\|^2_{L^2(\Om)}
\\
& \geq \, \al\Big(C_1 \|v\|_{W^{2,p}(\Om)}- C_1 \|g\|_{W^{2,p}(\Om)} -C_2\Big)\, +\frac{\beta}{2} \|\D^{\bar n} v\|^2_{L^2(\Om)}
\\
& \geq \, \al C_1 \|v\|_{W^{2,2}(\Om)} +\frac{\beta}{2} \|\D^{\bar n} v\|^2_{L^2(\Om)} - \al C_1 \|g\|_{W^{2,p}(\Om)} -\al C_2
\end{split}
\]
which implies
\[
\al C_1 \|v\|_{W^{2,2}(\Om)} +\frac{\beta}{2} \|\D^{\bar n} v\|^2_{L^2(\Om)} \, \leq \, \E_{\infty,\al,\be}(g) \,+\,  \al C_1 \|g\|_{W^{2,p}(\Om)} +\al C_2, 
\]
for any $v\in (W^{\bar n, 2}\cap W^{1,2}_g )(\Om)$. Now, by Poincar\'e inequality in $W^{1,2}_g(\Om)$ we have
\[
\|v\|_{L^2(\Om)} \, \leq\, C \Big(\|\D v\|_{L^2(\Om)} \,+\, \|g\|_{W^{1,2}(\Om)}\Big)
\]
for some $C>0$, and by the interpolation inequalities in the Sobolev space $W^{\bar n,2}(\Om)$, we have
\[
\|\D^k v\|_{L^2(\Om)} \, \leq\, C \Big(\|\D v\|_{L^2(\Om)} \,+\, \|\D^{\bar n} v\|_{L^2(\Om)}\Big),
\]
for some $C>0$ and any $k\in\{1,...,\bar n\}$. By putting the last three estimates together, we conclude that
\[
\| v\|_{W^{\bar n ,2}(\Om)} \, \leq \, C
\]
where the constant $C>0$ in general depends on $p$ but is uniform for $v \in (W^{\bar n, 2}\cap W^{1,2}_g )(\Om)$.

Let now $(u_p^m)_1^\infty$ be a minimising sequence of $\E_{p,\al,\be}$:
\[
\E_{p,\al,\be}(u_p^m) \larrow \inf\Big\{\E_p(v)\ : \ v\in (W^{\bar n, 2}\cap W^{1,2}_g )(\Om) \Big\},
\]
as $m\to \infty$. Then, by the coercivity estimate, we have the uniform bound
\[
\| u_p^m\|_{W^{\bar n ,2}(\Om)} \, \leq \, C
\]
for some $C>0$ independent of $m\in\N$. By standard weak and strong compactness arguments in Sobolev and H\"older spaces, together with the Morrey estimate
\[
\| v \|_{C^{k-[n/2]-1,\si}(\overline{\Om})}\, \leq\, C \| v\|_{W^{k ,2}(\Om)},
\]
applied to $k=\bar n$, there exists a subsequence $(u_p^{m_k})_1^\infty$ and $u_p\in (W^{\bar n, 2}\cap W^{1,2}_g )(\Om)$ such that, along this subsequence we have
\[
\left\{\ \ 
\begin{array}{ll}
u_p^m \larrow u_p, & \text{ in }C^{2}(\overline{\Om}),
\smallskip
\\
\D^k u_p^m \larrow \D^k u_p, & \text{ in }L^{2}\big(\Om,\R^{n^{\ot k}}_s\big), \text{ for all }k\in\{3,...,\bar n-1\},
\smallskip
\\
\D^{\bar n} u_p^m \weak \D^{\bar n} u_p, & \text{ in }L^{2}\big(\Om,\R^{n^{\ot \bar n}}_s\big),
\smallskip
\end{array}
\right.
\]
as $m_k\to \infty$. The above modes of convergence and the continuity of the function $K$ defining the operator $\K$ imply that $\K[u^m_p] \larrow \K[u_p]$ uniformly on $\mathcal{K}$ as  $m_k \to \infty$. Therefore,
\[
\big\||\K[u^m_p]-k^\ga|_{(p)} \big\|_{L^p(\mathcal{K},\mH^\kappa)} \larrow \big\||\K[u_p]-k^\ga|_{(p)} \big\|_{L^p(\mathcal{K},\mH^\kappa)}
\]
as $m_k \to \infty$. Additionally, by the continuity of the function $F$ defining the operator $\F$ and the uniform convergence of the minimising sequence up to second order derivatives, we have
\[
\F[u^m_p] \larrow \, \F[u_p]\ \ \text{ in }C(\overline{\Om}),
\]
as $m_k \to \infty$. Finally, by weak lower semi-continuity in $L^2$ we have
\[
\big\| \D^{\bar n} u_p   \big\|_{L^2(\Om)}\, \leq \, \liminf_{k\to\infty} \, \big\|\D^{\bar n} u^{m_k}_p \big\|_{L^2(\Om)}.
\]
By putting all the above together, we infer that
\[
\E_{p,\al,\be} (u_p) \, \leq \, \liminf_{k\to\infty} \, \E_{p,\al,\be}(u^{m_k}_p) \,\leq\, \inf\Big\{\E_{p,\al,\be}(v)\ : \ v\in (W^{\bar n, 2}\cap W^{1,2}_g)(\Om) \Big\},
\]
which concludes the proof.
\ep

Note that the proof above reveals the fact that $\E_{p,\al,\be}$ is weakly lower semi-continuous on the space $(W^{\bar n, 2}\cap W^{1,2}_g)(\Om)$, even though it is not explicitly stated.


\begin{lemma} \label{lemma2} For any $\al,\be,\ga>0$, there exists a (global) minimiser $u_\infty$ of $\E_{\infty,\al,\be}$ in the space $(W^{\bar n, 2}\cap W^{1,2}_g)(\Om)$, as well as a sequence of minimisers $(u_{p_i})_1^\infty$ of the respective $\E_{p,\al,\be}$-functionals from Lemma \ref{lemma1}, such that \eqref{1.17} holds true.
\end{lemma}

\bp For each $p>n$, let $u_p \in (W^{\bar n, 2}\cap W^{1,2}_g)(\Om)$ be the minimiser of $\E_{p,\al,\be}$ given by Lemma \ref{lemma1}. (We will follow a similar method and utilise the estimates appearing therein.) For any fixed $q \in (n,\infty)$ and $p\geq q$, H\"older's inequality and minimality yield
\[
\E_{q,\al,\be}(u_p)\, \leq\, \E_{p,\al,\be}(u_p)\, \leq\, \E_{p,\al,\be}(g)\, \leq\, \E_{\infty,\al,\be}(g)\, <\, \infty.
\]
By employing again the coercivity of $\E_{q,\al,\be}$, we have the estimate
\[
\begin{split}
\E_{q,\al,\be}(u_p) \, \geq\, \al \Big(C_1\|u_p\|_{W^{2,q}(\Om)}-C_1\|g\|_{W^{2,q}(\Om)}-C_2 \Big) + \frac{\be}{2} \|\D^{\bar n} u_p\|^2_{L^{2}(\Om)} ,
\end{split}
\]
which by Poincar\'e's inequality and the interpolation inequalities in $W^{\bar n,2}(\Om)$, yield
\[
\sup_{p\geq q} \, \| u_p \|_{W^{\bar n,2}(\Om)}\, \leq\, C
\]
for some $C>0$ depending on $q$ and all the parameters, but independent of $p$. By a standard diagonal argument, for any sequence $(p_i)_1^\infty$ with $p_i \larrow \infty$ as $i\to \infty$, there exists a function
\[
u_\infty \in (W^{\bar n, 2}\cap W^{1,2}_g)(\Om)
\]
and a subsequence (denoted again by $(p_i)_1^\infty$) along which \eqref{1.17} holds true. It remains to show that $u_\infty$ is in fact a minimiser of $\E_\infty$ over the same space. To this end, note that for any fixed $q\in (n,\infty)$ and $p\geq q$, we have
\[
\E_{q,\al,\be}(u_p)\, \leq\, \E_{p,\al,\be}(u_p)\, \leq\, \E_{p,\al,\be}(v)\, \leq\, \E_{\infty,\al,\be}(v)
\]
for any $v\in (W^{\bar n, 2}\cap W^{1,2}_g)(\Om)$. By the weak lower semi-continuity of $\E_{q,\al,\be}$ in the space $(W^{\bar n, 2}\cap W^{1,2}_g)(\Om)$ demonstrated in Lemma \ref{lemma1}, we have
\[
\E_{q,\al,\be}(u_\infty)\, \leq\, \liminf_{i\to \infty} \, \E_{q,\al,\be}(u_{p_i})\, \leq\, \E_{\infty,\al,\be}(v),
\]
for any $v\in (W^{\bar n, 2}\cap W^{1,2}_g)(\Om)$. 
By letting $q\to \infty$ in the estimate above, we deduce that 
\[
\E_{\infty,\al,\be}(u_\infty)\,\, \leq\, \inf\Big\{ \E_{\infty,\al\be}(v) \ : \ v \in (W^{\bar n, 2}\cap W^{1,2}_g)(\Om) \Big\},
\]
as desired.
\ep


\begin{lemma} \label{lemma3} For any $\al,\be,\ga>0$ and $p>n$, consider the minimiser $u_p$ of the functional $\E_{p,\al,\be}$ over $(W^{2,p}\cap W^{1,p}_g)(\Om)$ constructed in Lemma \ref{lemma1}. Consider also the signed Radon measures $\mu_p \in \mM(\Om)$ and $\nu_p \in \mM(\mathcal{K})$, defined in \eqref{1.18}. Then, the triplet $(u_p,\mu_p,\nu_p)$ satisfies the PDE \eqref{1.19} in the distributional sense, namely \eqref{1.21-A} holds true for all test functions $\phi \in C^{\bar n}_c(\Om)$. 
\end{lemma}

\bp We involve a standard Gateaux differentiability argument. Let us begin by checking that $\mu_p,\nu_p$ as defined in \eqref{1.18} are uniformly bounded Radon measures when $u_p \in W^{\bar n,2}(\Om)$. Since by the regularity of $F,K,u_p$ they obviously define absolutely continuous measures, it suffices to check that by H\"older inequality's, we have the total variation estimates
\[
\begin{split}
\|\nu_p\|(\mathcal{K})\, & \leq \, \Big(\big \| |\K[u_p]-k^\ga|_{(p)} \big\|_{L^p(\mathcal{K},\mH^\kappa)}\Big)^{1-p} \, {\av_\mathcal{K}} \big|\K[u_p]-k^\ga \big|^{p-1}_{(p)}  \, \mathrm{d} \mH^\kappa
\\
&\leq\, \Big(\big \| |\K[u_p]-k^\ga|_{(p)} \big\|_{L^p(\mathcal{K},\mH^\kappa)}\Big)^{1-p} \left(\, \, {\av_\mathcal{K}} \big|\K[u_p]-k^\ga \big|^{p}_{(p)}  \, \mathrm{d} \mH^\kappa \right)^{\!\!\frac{p-1}{p}}
\\
& = \, 1
\end{split}
\]
and similarly
\[
\begin{split}
\|\mu_p\|(\Om)\, & \leq \, \Big(\big \| | \F[u_p] |_{(p)} \big\|_{L^p(\Om)}\Big)^{1-p} \, {\av_\Om} \big| \F[u_p] \big|^{p-1}_{(p)}  \, \mathrm{d} \mL^n
\\
&\leq\, \Big(\big \| | \F[u_p] |_{(p)} \big\|_{L^p(\Om)}\Big)^{1-p} \left(\, \, {\av_\Om} \big|\F[u_p] \big|^{p}_{(p)}  \, \mathrm{d} \mL^n \right)^{\!\!\frac{p-1}{p}}
\\
& = \, 1.
\end{split}
\]
Next, fix $\phi \in C^{\bar n}_c(\Om)$. By using the regularity assumptions on $F,K$, we formally compute, recalling the abbreviations $\F[v]=F(\cdot,v,\D v, \D^2v)$ and $\K[v]=K(\cdot,v,\D v)$:
\[
\begin{split}
\frac{\mathrm{d}}{\mathrm{d}\e}\bigg|_{\e=0}   & \!\! \E_{p,\al,\be}(u_p+\e\phi) 
\\
&= 
\! \left(\, \, {\av_\mathcal{K}} \big|\K[u_p]-k^\ga \big|^{p}_{(p)}  \, \mathrm{d} \mH^\kappa \right)^{\!\!\frac{1}{p}-1} \!\! \, {\av_\mathcal{K}} \big|\K[u_p]-k^\ga \big|^{p-2}_{(p)}\big( \K[u_p]-k^\ga \big)   \centerdot
\\
& \ \ \  \centerdot \Big[ \K_{\mathrm r} [u_p] \,\phi \, +\, \K_{\mathrm p}[u_p]\cdot \D \phi  \Big] \, \mathrm{d} \mH^\kappa
\\
& \ \ \ +\, \al \left(\, \, {\av_\Om} \big| \F[u_p] \big|^{p}_{(p)}  \, \mathrm{d} \mL^n\right)^{\!\!\frac{1}{p}-1} \!\! \, {\av_\Om} \big| \F[u_p] \big|^{p-2}_{(p)} \, \F[u_p] \centerdot
\\
&  \ \ \ \centerdot \Big[\F_{\mathrm r} [u_p]\phi+\F_{\mathrm p} [u_p]\cdot\D \phi+\F_{\mathrm X} [u_p]:\D^2 \phi\Big] \, \mathrm{d} \mL^n
\\
& \ \ \ +\, \beta \, {\, {\av_\Om}} \, \D^{\bar n} u_p : \D^{\bar n} \phi \, \mathrm{d} \mL^n.
\end{split}
\]
By invoking the definitions of $\mu_p,\nu_p$, the above yields that
\[
\begin{split}
\frac{\mathrm{d}}{\mathrm{d}\e}\bigg|_{\e=0} \!\! \E_{p,\al,\be}(u_p+\e\phi) 
\, & = 
\int_{\mathcal K} \Big[ \K_{\mathrm r} [u_p] \,\phi \, +\, \K_{\mathrm p}[u_p]\cdot \D \phi  \Big] \, \mathrm{d} \nu_p
\\
& \ \ \ +\, \al \int_\Om \Big[\F_{\mathrm r} [u_p]\phi+\F_{\mathrm p} [u_p]\cdot\D \phi+\F_{\mathrm X} [u_p]:\D^2 \phi\Big] \, \mathrm{d} \mu_p
\\
& \ \ \ +\, \beta \, {\, {\av_\Om}} \, \D^{\bar n} u_p : \D^{\bar n} \phi \, \mathrm{d} \mL^n.
\end{split}
\]
Since $u_p$ is the minimiser of $\E_{p,\al,\be}$, we have that $\E_{p,\al,\be}(u_p)\leq \E_{p,\al,\be}(u_p+\e \phi)$ for all $\e\in\R$ and all $\phi \in C^{\bar n}_c(\Om)$. Hence, our computation implies that the PDE \eqref{1.19} is indeed satisfied as claimed, once we confirm that the formal computation is rigorous, and that therefore $\E_p$ is Gateaux differentiable at the minimiser $u_p$ for any direction $\phi \in C^{\bar n}_c(\Om)$ because by the continuity of $\F,\K$ and the fact that $u_p \in(C^2 \cap W^{\bar n,2})(\Om)$, $\F[u_p] \in C(\overline{\Om})$ and $\K[u_p]-k^\ga \in L^\infty(\mathcal{K},\mH^\kappa)$, H\"older's inequality implies that 
\[
\F_{\mathrm r} [u_p]\phi+\F_{\mathrm p} [u_p]\cdot\D \phi+\F_{\mathrm X} [u_p]:\D^2 \phi \ \in \, C(\overline{\Om})
\]
and 
\[
 K_{\mathrm r} (\cdot,u_p,\D u_p) \,\phi + K_p(\cdot,u_p,\D u_p)\cdot \D \phi  \ \in \ C(\mathcal{K}),
\]
for any $\phi \in C^{\bar n}_c(\Om)$. Finally, $\D^{\bar n} u_p : \D^{\bar n} \phi  \in L^1(\Om)$ since $\D^{\bar n} u_p \in L^2\big(\Om,\R^{n^{\ot \bar n}}_s\big)$.
\ep

\begin{lemma} \label{lemma4} For any $\al,\be,\ga>0$, consider the minimiser $u_\infty$ of $\E_{\infty,\al,\be}$ constructed in Lemma \ref{lemma2} and the minimisers $(u_p)_{p>n}$ of the functionals $(\E_{p,\al,\be})_{p>n}$. Then, there exist signed Radon measures $\mu_\infty \in \mM(\Om)$ and $\nu_\infty \in \mM(\mathcal{K})$ such that the triplet $(u_\infty,\mu_\infty,\nu_\infty)$ satisfies the PDE \eqref{1.14} in the distributional sense, that is \eqref{1.16-A} holds true. Additionally, there exists a further subsequence along which the weak* modes of convergence of \eqref{1.18} hold true as $p\to \infty$.
\end{lemma}

\bp By the proof of Lemma \ref{lemma3}, we have the uniform in $p$ total variation bounds $\|\mu_p\|(\Om)\leq 1$ and  $\|\nu_p\|(\mathcal{K})\leq 1$. Hence, by the sequential weak* compactness of the spaces of Radon measures
\beq \label{duality}
\mM(\Om)\, =  \big(C_0(\Om)\big)^*, \ \ \ \mM(\mathcal K)\, =  \big(C(\mathcal{K})\big)^*,
\eeq
there exists a further subsequence denoted again by $(p_i)_1^\infty$ such that $\mu_p \weakstar \mu_\infty$ in $\mM(\Om)$ and  $\nu_p \weakstar \nu_\infty$ in $\mM(\mathcal{K})$, as $p_i \to \infty$. Fix now $\phi \in C^{\bar n}_c(\Om)$. By Lemma \ref{lemma3}, we have that the triplet $(u_p,\mu_p,\nu_p)$ satisfies \eqref{1.19}, that is \eqref{1.21-A} holds true for any fixed test function. Since $u_p \larrow u_\infty$ in $C^2(\overline{\Om})$ as $p_i \to \infty$, $F,F_{\mathrm r},F_{\mathrm p}, F_{\mathrm X},K_{\mathrm r}, K_{\mathrm p}$ are all continuous up to the boundary, we have as $p_i \to \infty$ that
\[
\begin{split}
 \K_{\mathrm r} [u_p]\,\phi + \K_{\mathrm p}[u_p] \cdot \D \phi \, & \larrow \, \K_{\mathrm r} \K_{\mathrm r} [u_\infty]\,\phi  + \K_{\mathrm p}[u_\infty] \cdot \D \phi,
\\
 \F_{\mathrm r} [u_p]\,\phi + \F_{\mathrm p}[u_p] \cdot \D \phi+ \F_{\mathrm X}[u_p] : \D^2 \phi \, &\larrow \, \F_{\mathrm r} [u_\infty]\,\phi  + \F_{\mathrm p}[u_\infty] \cdot \D \phi \F_{\mathrm X}[u_\infty] : \D^2 \phi,
\end{split}
\]
uniformly on $\mathcal K$ and on $\overline{\Om}$ respectively, for any fixed $\phi \in C^{\bar n}_c(\Om)$. Further, we have $\mu_p\weakstar \mu_\infty$ and $\nu_p\weakstar \nu_\infty$ in $\mM(\Om)$ and in $\mM(\mathcal K)$ respectively as $p_j\to \infty$. By standard properties of the weak*-strong continuity of the duality pairings \eqref{duality}, we have that the first two terms of \eqref{1.16-A} converge to the respective first two terms of \eqref{1.21-A}, as $p_i \to \infty$. Finally, since $\D^{\bar n}u_p \weak \D^{\bar n}u_p$ in $L^2\big(\Om,\R^{n^{\ot \bar n}}_s\big)$ and $ \D^{\bar n}\phi \in \smash{L^2\big(\Om,\R^{n^{\ot \bar n}}_s\big)}$, we also have that the last term of \eqref{1.16-A} converges to the respective last term of \eqref{1.21-A}, as $p_i \to \infty$. Hence, we have indeed obtain \eqref{1.14} by passing to the limit as  $p_i \to \infty$ in \eqref{1.19}.
\ep

\begin{lemma} \label{lemma5} For any $\al,\be,\ga>0$, $p>n$ and $u^0 \in (W^{2,p}\cap W^{1,p}_g)(\Om)$ such that
\[
\big\| k^\ga - \K[u^0] \big\|_{L^\infty(\mathcal{K},\mH^\kappa)} \,\leq\, \ga,
\]
the ($(\al,\be,\ga)$-dependent) minimiser $u_p$ of $\E_p$ (constructed in Lemmas \ref{lemma1}-\ref{lemma4}), satisfies the error bounds \eqref{1.24}. If additionally $u^0 \in (W^{\bar n,2}\cap W^{1,2}_g)(\Om)$, then the ($(\al,\be,\ga)$-dependent) minimiser $u_\infty$ of $\E_\infty$ (constructed in Lemmas \ref{lemma1}-\ref{lemma4}), satisfies the error bounds \eqref{1.20}.
\end{lemma}

\bp Let us denote $k^0:= \K[u^0]$, noting that $k^0 \in C(\mathcal{K})$ and that
\[
\|k^\ga -k^0\|_{L^\infty(\mathcal{K},\mH^\kappa)} \, \leq \, \ga.
\]
Recall that for any $p>n$, $u_p$ is a global minimiser of $\E_{p,\al,\be}$ in $(W^{\bar n,2}\cap W^{1,2}_g)(\Om)$. Therefore,
\[
\E_{p,\al,\be}(u_p) \, \leq\, \E_{p,\al,\be}(u^0).
\]
This implies  
\[
\begin{split}
\big \|\K[u_p]-k^\ga \big\|_{L^p(\mathcal{K},\mH^\kappa)} & +\, \al \big \| \F[u_p] \big\|_{L^p( \Om)} + \, \frac{\be}{2} \, \|\D^{\bar n} u_p\|^2_{L^2(\Om)}
\\
&\leq\, \big \|\K[u^0]-k^\ga \big\|_{L^p(\mathcal{K},\mH^\kappa)}+ \al \big \| \F[u^0] \big\|_{L^p( \Om)} + \, \frac{\be}{2} \, \|\D^{\bar n} u^0\|^2_{L^2(\Om)}. 
\end{split}
\]
The above together with Minkowski's and H\"older's inequalities yields
\[
\begin{split}
\big \|\K[u_p]- \K[u^0] \big\|_{L^p(\mathcal{K},\mH^\kappa)} & \leq\, \big \|\K[u^0]-k^\ga \big\|_{L^p(\mathcal{K},\mH^\kappa)}+\, \big \|\K[u^0]-k^\ga \big\|_{L^p(\mathcal{K},\mH^\kappa)}\,
\\
& \ \ \ \ +\, \al \big \| \F[u^0] \big\|_{L^p( \Om)} + \, \frac{\be}{2} \, \|\D^{\bar n} u^0\|^2_{L^2(\Om)}
\\
& = \, 2 \|k^\ga -k^0\|_{L^p(\mathcal{K},\mH^\kappa)}   + \al \big \| \F[u^0] \big\|_{L^p(\Om)} + \, \frac{\be}{2} \, \|\D^{\bar n} u^0\|^2_{L^2(\Om)}
\\
& \leq 2\ga + \al \, \| \F[u^0]\|_{L^p(\Om)} + \, \frac{\be}{2} \, \|\D^{\bar n} u^0\|^2_{L^2(\Om)},
\end{split}
\]
as claimed. To obtain the corresponding estimate for $u_\infty$ in the case that additionally $u^0 \in (W^{\bar n,2}\cap W^{1,2}_g)(\Om)$, we pass to the limit as $p_i\to \infty$ in the estimate above and the conclusion follows by letting $p_i\to \infty$ due to the strong convergence $u_p \larrow u_\infty$ in $C^2(\overline{\Om})$ as $p_i\to \infty$. The lemma ensues.
\ep

The proof is now complete by noting that the statements \eqref{1.15}-\eqref{1.16} in Theorem \ref{theorem1} follow from Proposition \ref{prop} and the established modes of convergence.

%
%
%
%
%
%
%
%
%
%
%
%
%
%


\bibliographystyle{amsplain}

\end{document}